\documentclass[12pt]{article}
\usepackage{amstex}
\usepackage{amssymb}
\usepackage{amscd}

\def\no{\noindent}
\newcommand{\qed}{{\unskip\nobreak\hfil
        \penalty50\hskip1em\hbox{}\nobreak\hfil
        $\square$\parfillskip=0pt\finalhyphendemerits=0 \par}}

\newtheorem{dfn}{Definition}[section]

\newtheorem{thm}[dfn]{Theorem}

\newtheorem{lem}[dfn]{Lemma}

\newtheorem{sublem}[dfn]{Sublemma}
\newtheorem{prop}[dfn]{Proposition}
\newtheorem{cor}[dfn]{Corollary}

\newtheorem{fact}[dfn]{Fact}

\def\proof{\par\medskip\noindent{\it Proof. }}
\def\lra{\longrightarrow}
\def\Lra{\Longrightarrow}
\def\ra{\rightarrow}

\def\0{\emptyset}
\def\R{{\Bbb R}}
\def\E{{\Bbb E}}
\def\H{{\Bbb H}}
\def\Z{{\Bbb Z}}

\def\N{{\Bbb N}}

\def\eps{\epsilon}
\def\al{\alpha}

\def\ga{\gamma}
\def\Ga{\Gamma}
\def\de{\delta}

\def\8{\infty}
\def\<{\langle}
\def\>{\rangle}

\def\BE{\begin{equation}}
\def\EE{\end{equation}}

\def\qi{quasi-isometry\ }

\def\geo{\partial_{\infty}}
\def\tits{\partial_{Tits}}

\def\ol{\overline}
\newcommand{\restr}{\mbox{\Large \(|\)\normalsize}}
\def\trip{\partial^3}

\begin{document}
\title{Groups quasi-isometric to symmetric spaces}
\author{
\setcounter{footnote}{1}
Bruce Kleiner\thanks{Partially
supported by NSF grants
DMS-95-05175 and DMS-96-26911.} 
and 
\setcounter{footnote}{0}
Bernhard Leeb\thanks{
Supported by SFB 256 (Bonn).}
}
\date{October 3,1996}
\maketitle

\begin{abstract}
\no
We determine the structure of finitely generated groups 
which are quasi-isometric to symmetric spaces of noncompact type,
allowing Euclidean de Rham factors
If $X$ is a symmetric space of noncompact type with no Euclidean
de Rham factor, and $\Ga$ is a finitely generated group quasi-isometric
to the product $\E^k\times X$, then there is an exact sequence
$1\ra H\ra\Ga\ra L\ra 1$ where $H$ contains a finite index copy
of $\Z^k$ and $L$ is a uniform lattice in the isometry group of
$X$.
\end{abstract}

\section{Introduction}

The main result of this paper is the following theorem.

\begin{thm}
\label{mainthm}
Let $X$ be a symmetric space of noncompact type with no Euclidean 
de Rham factor,
and let $Nil$ be a simply connected nilpotent Lie group 
equipped with a left-invariant Riemannian metric. 
Suppose that $\Ga$ is a finitely generated group quasi-isometric to 
$Nil\times X$. 
Then there is an exact sequence
\begin{equation}
\label{exseq}
1\lra H\lra \Ga\lra L\lra 1
\end{equation}
where $H$ is a finitely generated group quasi-isometric to $Nil$ 
and $L$ is a uniform lattice in the isometry group of $X$. 
\end{thm}
In particular, 
when $Nil$ is the trivial group then $\Ga$
is a finite extension of a uniform lattice
in $Isom(X)$, and when $Nil\simeq \R^k$ then 
$H$ is virtually abelian of rank $k$ by \cite{polygro,Pansu1}. 

Certain cases of theorem \ref{mainthm}
were proved earlier.  \cite{Tukia,hypmfds} determined the 
structure of groups quasi-isometric to $\H^{n\geq 3}$,
\cite{Pansu} handled the other rank 1 cases besides $\H^2$,
and \cite{CaJu,Gab} settled the $\H^2$ case.  \cite{Rieffel}
studied groups quasi-isometric to $\H^2\times\R$; in this
case it was already known \cite{Gersten}
that such groups need not admit  discrete cocompact isometric
actions on
$\H^2\times\R$.

It is an intriguing problem to classify the finitely
generated groups
quasi-isometric to a given space/group. 
Some other cases where the classification is known
are: free groups \cite{hypgps,GhysHar},  certain
nilpotent groups \cite{polygro,Pansu1},
non-uniform
lattices in rank $1$ symmetric spaces \cite{Schwartz}
and in symmetric spaces without rank $1$ factors
\cite{Eskin}, and fundamental groups of Haken 3-manifolds
with nontrivial geometric decomposition \cite{KapLe}.

\bigskip
We gratefully acknowledge support by 
the RiP-program at the Mathematisches Forschungsinstitut Oberwolfach. 

\tableofcontents

\section{Preliminaries}

In this section we recall some basic definitions
and notation.  See
\cite{AsyInv} for more discussion and background.

\begin{dfn}
A map $f:X\lra Y$ between metric spaces is an
{\bf $(L,A)$ quasi-isometry} if for every $x_1,x_2\in X$
$$L^{-1}d(x_1,x_2)+A\leq d(x_1,x_2)\leq Ld(x_1,x_2)+A,$$
and for every $y\in Y$ we have $d(y,f(X))<A$.  Two 
quasi-isometries $f_1,f_2:X\lra Y$ are {\bf equivalent}
if $d(f_1,f_2)<\infty$.
\end{dfn}

If $\Ga$ is a finitely generated group, then any two 
word metrics on $\Ga$ are biLipschitz to one another
by $id_\Ga:\Ga\ra \Ga$.  We will implicitly endow
our finitely generated groups with word metrics.

\begin{dfn}
An {\bf $(L,A)$-quasi-action} of a group $\Ga$ on a metric
space $Z$ is a map $\rho:\Ga\times Z\ra Z$ 
so that $\rho(\ga,\cdot):Z\ra Z$ is an $(L,A)$
quasi-isometry for every $\ga\in \Ga$,
$d(\rho(\ga_1,\rho(\ga_2,z)),\rho(\ga_1\ga_2,z))<A$ for 
every $\ga_1,\ga_2\in\Ga$, $z\in Z$, and $d(\rho(e,z),z)<A$
for every $z\in Z$.
\end{dfn}

We will denote the self-map $\rho(\ga,\cdot):Z\ra Z$ by $\rho(\ga)$. 
$\rho$ is {\bf discrete} if for any point $z\in Z$ 
and any radius $R>0$, the set of all $\ga\in\Ga$ such that $\rho(\ga,z)$ is contained 
in the ball $B_R(z)$ is finite. 
$\rho$ is {\bf cobounded} if $Z$ coincides with a 
finite tubular neighborhood of the ``orbit'' $\rho(\Ga) z\subset Z$ 
for every $z$.  
If $\rho$ is a discrete cobounded quasi-action 
of a finitely generated group $\Ga$ on a geodesic metric space $Z$, it
follows easily that 
the map $\Ga\ra Z$ given by $\ga\mapsto\rho(\ga,z)$ 
is a quasi-isometry for every $z\in Z$.

\begin{dfn}
Two quasi-actions $\rho$ and $\rho'$ are {\bf equivalent} 
if there exists a constant $D$ so that 
$d(\rho(\ga),\rho'(\ga)<D$ for all $\ga\in\Ga$. 
\end{dfn}

\begin{dfn}
Let $\rho$ and $\rho'$  be a quasi-actions of $\Ga$ on $Z$
and $Z'$ respectively, and let $\phi:Z\ra Z'$ be a \qi. 
Then $\rho$ is {\bf quasi-isometrically conjugate
to $\rho'$ via  $\phi$}
if there is a $D$ so that 
$d(\phi\circ\rho(\ga),\rho'(\ga)\circ\phi)<D$
for all $\ga\in\Ga$.
\end{dfn}

\begin{lem}
(cf \cite[8.2.K]{hypgps})
\label{hyplem}
Let $X$ be a Hadamard manifold of dimension $\geq 2$
with sectional curvature $\leq K<0$, and let $\geo X$ 
denote the geometric boundary of $X$ with the cone
topology.  Recall that every quasi-isometry 
$\Phi:X\lra X$ induces a boundary homeomorphism
$\geo\Phi:\geo X\ra\geo X$. 
\begin{enumerate}

\item
\label{discretecobounded}
 If $\rho:\Ga\times X\ra X$ is a quasi-action on $X$,
then $\rho$ is discrete (respectively cobounded) iff
$\geo\phi$ acts properly discontinuously (respectively
cocompactly) on the space of distinct triples in $\geo X$.

\item
\label{boundarymaps}
 Given $(L,A)$ there is a $D$ so that if $\phi_k$, $\psi$ are $(L,A)$
quasi-isometries, then $\geo\phi_k$ converges uniformly to $\geo\psi$ iff
 $\limsup d(\phi_k x,\psi x)<D$ for 
every $x\in X$. 
In particular, if $\phi_1,\phi_2:X\lra X$ are $(L,A)$ quasi-isometries
with the same boundary mappings,
then $d(\phi_1,\phi_2)<D$.
\end{enumerate}
\end{lem}

\proof 
Let $\trip X\subset\geo X\times\geo X\times\geo X$
denote the subspace of distinct triples.  The uniform negative
curvature of $X$ implies that there is a $D_0$ 
depending only on $K$ such that\\
(a) For every $x\in X$ there is a triple $(\xi_1,\xi_2,\xi_3)
\in\trip X$ such that $d(x,\ol{\xi_i\xi_j})<D_0$ for every
$1\leq i\neq j\leq 3$, where $\ol{\xi_i\xi_j}$
denotes the geodesic with ideal endpoints $\xi_i$, $\xi_j$.
Moreover for every $C$ the set $\{(\xi_1,\xi_2,\xi_3)
\mid \mbox{$d(x,\ol{\xi_i\xi_j})<C$ for all $1\leq i\neq j\leq 3$}\}$
has compact closure in $\trip X$.\\
and\\
(b) For every $(\xi_1,\xi_2,\xi_3)\in\trip X$
there is a point $x\in X$ so that $d(x,\ol{\xi_i\xi_j})<D_0$
for each $1\leq i\neq j\leq 3$.  And for every $C$
there is a $C'$ depending only on $C$ and $K$
so that $\{x\in X\mid\mbox{$d(x,\ol{\xi_i\xi_j})<C$ for every
$1\leq i\neq j\leq 3$}\}$ has diameter $<C'$.

\smallskip
\no
\ref{discretecobounded} and \ref{boundarymaps} follow
easily from this.
\qed

\section{Projecting quasi-actions to the factors}
\label{projectingqactions}

Let $Nil$ and $X$ be as in Theorem \ref{mainthm} 
and decompose $X$ into irreducible factors: 
\begin{equation}
X=\prod_{i=1}^l X_i
\end{equation}
Suppose $\rho$ is a quasi-action of the finitely generated group $\Ga$ 
on $Nil\times X$. 
We denote by $p:Nil\times X\ra X$ the canonical projection.
By applying \cite[Theorem 1.1.2]{KlLe} to each quasi-isometry 
$\rho(\ga)$ we construct 
quasi-actions $\rho_i$ of $\Ga$ on $X_i$ 
so that 
\[
d(p\circ\rho(\ga),\prod_{i=1}^k\rho_i(\ga)\circ p)<D
\]
for all $\ga\in\Ga$ and some positive constant $D$.

\section{Straightening cocompact quasi-actions on irreducible
symmetric spaces}

\bigskip
The following result is a direct consequence of 
\cite[Th\'eor\`eme 1]{Pansu} and \cite[Theorem 1.1.3]{KlLe}. 

\begin{fact}
\label{qiconjrigidcase}
Let $X$ be an irreducible symmetric space 
other than a real or complex hyperbolic space. 
Then every quasi-action on $X$ is equivalent to an isometric action.
\end{fact}
\proof
Let $\rho$ be a quasi-action of a
group $\Ga$ on $X$.  By the results just cited, 
there is an isometry $\bar\rho(\ga)$ at finite distance 
from the  quasi-isometry $\rho(\ga)$ for every $\ga\in \Ga$. 
This isometry is unique and its distance from $\rho(\ga)$ is 
uniformly bounded\footnote{The uniformity in the rank one case 
follows lemma \ref{hyplem}.}
in terms of the constants of the quasi-action. So $\bar\rho$
is an isometric action equivalent to $\rho$.
$\Box$

We recall that the real and complex hyperbolic spaces of all dimensions 
admit quasi-isometries which are not equivalent to isometries
\cite{Pansu}. 

\begin{fact}
\label{qiconjrealcomplexhyp}
Any cobounded quasi-action $\rho$ on a real or complex hyperbolic space 
is quasi-isometrically conjugate to an isometric action. 
\end{fact}

This result is proven in \cite{Tukia} in the real-hyperbolic case. 
Using Pansu's theory of Carnot differentiability one can carry out 
Tukia's arguments for all rank-one symmetric spaces 
other than hyperbolic plane, cf.\ \cite[sec.\ 11]{Pansu}. 
Another proof for the complex-hyperbolic case can be found in \cite{Chow}.

\begin{fact}
\label{qiconjhypplane}
Let $\rho$ be a cobounded quasi-action of a group $\Ga$
on $\H^2$.  Then $\rho$ is quasi-isometrically conjugate
to a cocompact isometric action of $\Ga$ on $\H^2$.
\end{fact}
\proof 
We recall that every \qi $\phi:\H^2\ra \H^2$
induces a quasi-symmetric homeomorphism $\geo\phi:\geo\H^2\ra\geo\H^2$,
see \cite{TukVai}; moreover the quasi-symmetry constant of $\geo\phi$
can be estimated in terms of the quasi-isometry constants of $\phi$.
Since equivalent quasi-isometries yield the same boundary
homeomorphism, every quasi-action $\rho$ on $\H^2$
induces a genuine action $\geo\rho$ on $\geo\H^2$ by uniformly 
quasi-symmetric homeomorphisms.

Let $\bar\Ga$ be the quotient of $\Ga$ by the kernel of the action $\geo\rho$,
and let $\pi:\Ga\ra\bar\Ga$ be the canonical epimorphism. 
If two elements $\ga_1,\ga_2\in\Ga$ have the same boundary map 
then  $d(\rho(\ga_1),\rho(\ga_2))$ is uniformly bounded by
lemma \ref{hyplem}.
Hence we may obtain a quasi-action $\bar\rho$ of $\bar\Ga$ on $\H^2$
by choosing $\ga\in\pi^{-1}(\bar\ga)$ for each $\bar\ga\in\bar\Ga$,
and setting $\bar\rho(\bar\ga)=\rho(\ga)$.  If $\bar\tau$ 
is an isometric action of $\bar\Ga$ on $\H^2$ and 
$\phi:\H^2\ra\H^2$ quasi-isometrically conjugates $\bar\rho$
into $\bar\tau$, then $\phi$ will quasi-isometrically conjugate
$\rho$ into the isometric action $\tau:\Ga\times\H^2\ra\H^2$
given by $\tau(\ga)=\bar\tau(\pi(\ga))$.  Hence it suffices
to treat the case when $\bar\Ga=\Ga$, and so we will assume that
$\geo\rho$ is an effective action.

\begin{lem}
The quasi-action $\rho$ is discrete if and only if 
the action $\geo\rho$ on $\geo\H^2$ is discrete in the compact-open topology. 
\end{lem}

\proof
Suppose $\geo\rho$ is discrete, and
let $(\ga_i)$ be a sequence in $\Ga$ so that $\rho(\ga_i)$ 
maps a point $p\in\H^2$ into a fixed ball $B_R(p)$.   Then
by a selection argument we may assume -- after passing
to a subsequence if necessary -- that there is a quasi-isometry
$\phi:\H^2\ra\H^2$ so that for every $q\in \H^2$
we have $\limsup_i d(\rho(\ga_i)(q),\phi(q))<D$ for some $D$.  Hence
the boundary maps $\geo\rho(\ga_i)$ converge to $\geo\phi$, 
and so the sequence $\geo\rho(\ga_i)$ is eventually constant.
Since $\rho$ is effective we conclude that $\ga_i$ is 
eventually constant.  Therefore $\rho$ is a discrete quasi-action.

If $\rho$ is a discrete quasi-action on $\H^2$, then 
$\geo\rho$ is discrete by lemma \ref{hyplem}.
\qed

\medskip
\no
{\em Proof of \ref{qiconjhypplane} continued.} 

\medskip
\no
{\em Case 1: $\geo\rho$ is discrete.}  
In this case, $\rho$ is a discrete convergence group action 
(lemma \ref{hyplem}) and by the work of \cite{CaJu,Gab}, there is
a discrete isometric action $\tau$ of $\Ga$
on $\H^2$ so that $\geo\rho$ is topologically conjugate
to $\geo\tau$.  Since $\rho$ is cobounded,
$\geo\rho$ acts cocompactly on the set of distinct triples of
points
in $\geo\H^2$ (lemma \ref{hyplem}); therefore
 $\geo\tau$ also acts cocompactly on the 
space of triples and so $\tau$ is a discrete, cocompact, isometric action
of $\Ga$ on $\H^2$.  
We now have two discrete, cobounded,
quasi-actions of $\Ga$ on $\H^2$, so they are quasi-isometrically
conjugate by some quasi-isometry $\psi:\H^2\ra\H^2$.

\medskip
\no
{\em Case 2: $\geo\rho$ is nondiscrete.} By \cite[Theorem 4]{Hin},
$\geo\rho$ is quasi-symmetrically conjugate to $\geo\tau$,
where $\tau$ is  an isometric action on $\H^2$.
  The conjugating quasi-symmetric
homeomorphism is the boundary of a quasi-isometry $\psi:\H^2\ra
\H^2$, \cite{TukVai}, which quasi-isometrically conjugates
$\geo\rho$ into the  isometric action action $\tau$.
Applying lemma \ref{hyplem} again, we conclude that
$\tau$ is cocompact.
\qed

Section \ref{projectingqactions}, 
\ref{qiconjrigidcase}, \ref{qiconjrealcomplexhyp} and \ref{qiconjhypplane}
imply: 

\begin{cor}
\label{conjqactions}
Let $X$ be a symmetric space of noncompact type without Euclidean factor. 
Then any cobounded quasi-action on $X$ is quasi-isometrically conjugate 
to a cocompact isometric action on $X$. 
\end{cor}

\section{A Growth estimate for small elements in nondiscrete 
cocompact subgroups
of $Isom(X)$}

\subsection{Parabolic isometries of symmetric spaces}
\label{paraboliccrap}

Let $X$ be a symmetric space of noncompact type, and let $G=Isom(X)$. 

\bigskip
An isometry $g\in G$ is {\bf semisimple} if its displacement function
$\de_g$ attains its infimum and {\bf parabolic} 
otherwise. 
 
\begin{lem}
Let $A\subset G$ be a finitely generated abelian group all of whose nontrivial elements 
are parabolic. 
Then $A$ has a fixed point at infinity.
\end{lem}
\proof
Recall that the nearest point projection to a closed convex subset is 
well-defined and distance non-increasing. 
This implies that if $C$ is a non-empty $A$-invariant closed convex set,
then for all displacement functions $\de_a$, $a\in A$, 
we have $\inf\de_a=\inf\de_a\restr_C$. 
Hence for all $n\in\N$, the intersection of the sublevel sets
$\{p\mid \de_{a_i}(p)\leq \inf\de_{a_i}+1/n\}$ is non-empty and contains a point $p_n$. 
We have $\de_{a_i}(p_n)\ra\inf\de_{a_i}$ 
for all $a_i$, and 
since the isometries $a_i$ are parabolic 
the sequence $\{p_n\}$ subconverges to an ideal boundary point 
$\xi\in\tits X$. 
It follows that the $a_i$ fix $\xi$. 
$\Box$

\begin{lem}
\label{kirsch}
Let $a_1,\dots,a_k\in Isom(X)$ be commuting parabolic isometries. 
Then there is a sequence of isometries $\{g_n\}\subset G$ 
so that for every $i$ the sequence $g_na_ig_n^{-1}$
subconverges to a semisimple isometry $\bar a_i$. 
\end{lem}
\proof
>From the proof of the previous lemma,
there is a sequence of points $\{p_n\}\subset X$ converging to an ideal point $\xi$ 
so that $\de_{a_i}(p_n)\ra\inf\de_{a_i}$ for all $a_i$. 
Pick isometries $g_n\in G$ such that $g_n\cdot p_n=p_0$. 
The conjugates $g_na_ig_n^{-1}$ have the same infimum displacement as $a_i$. 
Since 
\[ \de_{g_na_ig_n^{-1}}(p_0)= \de_{a_i}(p_n)\ra\inf\de_{a_i}\quad , \]
 the $g_na_ig_n^{-1}$ subconverge to a semisimple isometry. 
$\Box$

We call the isometry $g\neq e$ {\bf purely parabolic} if the identity 
is the only semisimple element in  
$\ol{Ad_G(G)\cdot g}$.

\subsection{The growth estimate}

\begin{prop}
\label{growth}
Let $X$ be a symmetric space of noncompact type with no Euclidean
de Rham factors.  Let $\Ga\subset G=Isom(X)$ be a finitely 
generated, nondiscrete, cocompact subgroup. Let $U\subset Isom(X)$
be a neighborhood of the identity, and
set $$f(k):=\#\{ g\in\Ga : \mbox{$|g|_\Ga < k$, $g\in U$}\},$$
where $|\cdot |_\Ga$ denotes a word norm on  $\Ga$.
Then $f$ grows faster than any polynomial, i.e. for every $d>0$
$\limsup_{k\ra\infty}\frac{f(k)}{k^d}=\infty$. 
\end{prop}
\proof
Let $\bar\Ga^o$ denote the identity component of the closure
of $\Ga$ in $G$.

\no
{\em Case 1: $\bar\Ga^o$ is nilpotent.} 
Let $A$ be the last non-trivial subgroup in the derived series of $\bar\Ga^o$. 
Then $A\subset\bar\Ga$ is a connected abelian subgroup of positive dimension,
$A$ is normal in $\bar\Ga$, 
and $\Ga\cap A$ is dense in $A$. 

\begin{lem}
For every $\de\in(0,1)$ there is a $\ga\in\Ga$ such that all eigenvalues 
of the automorphism 
$Ad_G(\ga)\restr_A:A\ra A$ have absolute value $<\de$. 
\end{lem}

\proof
See section \ref{paraboliccrap} for terminology.

\no
{\em Step 1: $A$ contains no semisimple isometries other
than $e$.}
Otherwise we can consider the intersection $C$ of the minimum sets 
for the displacement functions $\de_a$ where $a$ runs through all semisimple elements in $A$. 
$C$ is a nonempty convex subset of $X$ which splits metrically as
$C\cong\E^k\times Y$. 
The flats $\E^k\times\{y\}$ are the minimal flats preserved by all semisimple 
elements in $A$. 
Since $\Ga$ normalises $A$ it follows that $C$ is $\Ga$-invariant.
The cocompactness of $\Ga$ implies that $C=X$ and $k=0$
because $X$ has no Euclidean factor. 
This means that the semisimple elements in $A$ fix all points,
a contradiction.

\no
{\em  Step 2: All non-trivial isometries in $A$ are purely parabolic.} 
If $a\in A$, $a\neq e$, is not purely parabolic 
then there is a sequence of isometries $g_n$ so that 
$g_nag_n^{-1}$ converges to a semisimple isometry $\bar a\neq e$. 
We can uniformly approximate the $g_n$ by elements in $\Ga$, 
i.e.\ there exist $\ga_n\in\Ga$ and a bounded sequence $k_n\in G$ 
subconverging to $k\in G$ so that 
$\ga_n=k_ng_n$. 
Then $\ga_na\ga_n^{-1}=k_ng_nag_n^{-1}k_n^{-1}$ subconverges 
to the non-trivial semisimple element $k\bar ak^{-1}$. 
This contradicts step 1.

\no
{\em  Step 3:}
Pick a basis $\{a_1,\dots,a_k\}$ for  $A\simeq\R^k$. 
By lemma \ref{kirsch} there exist elements $g_n\in G$ so that 
$g_na_ig_n^{-1}\ra e$ for all $a_i$. 
We approximate the $g_n$ as above by $\ga_n$ so that the sequence 
$\ga_ng_n^{-1}$ is bounded. 
Then $\ga_na_i\ga_n^{-1}\ra e$ for all $a_i$. 
The lemma follows by setting $\ga=\ga_n$ for sufficiently large $n$.
$\Box$

\medskip
\no
{\em Proof of case 1 continued.} By the lemma, 
there is a $\ga\in\Ga$, $\ga\neq e$,  
and a norm $\|\cdot\|_A$ on $A$ such that for all $a\in A$ we have
\[ \|\ga a\ga^{-1}\|_A < \frac{1}{2} \|a\|_A   .\]
Consider a neighborhood $U$ of $e$ in $G$. 
Let $r>0$ be small enough so that 
$\{a\in A : \|a\|_A<r\}\subset U$ 
and pick $\al\in\Ga\cap A$ with $\|\al\|_A<r/2$. 
Then the elements 
\[   \ga_{\eps_0\dots\eps_{n-1}} = 
   \al^{\eps_0}\cdot(\ga\al\ga^{-1})^{\eps_1}\cdot\dots\cdot
   (\ga^{n-1}\al\ga^{1-n})^{\eps_{n-1}}      \]
for $\eps_i\in\{0,1\}$ are $2^n$ pairwise distinct elements contained in $\Ga\cap U$ 
with word norm 
$|\ga_{\eps_0\dots\eps_{n-1}}|_{\Ga}<n^2(|\al|_{\Ga}+|\ga|_{\Ga})$. 
This implies superpolynomial growth of $f$. 

\medskip\no
{\em Case 2: $\bar\Ga^o$ is not nilpotent.}
Define an increasing sequence (the upper central series) of nilpotent 
Lie subgroups $Z_i\subset\bar\Ga^o$ inductively as follows: 
Set $Z_0=\{e\}$ and let $Z_{i+1}$ be the inverse image in $\bar\Ga^o$ of the center 
in $\bar\Ga^o/Z_i$. 
The dimension of $Z_i$ stabilizes and we choose $k$ 
so that $\dim Z_k$ is maximal. 
Then the center of $\bar\Ga/Z_k$ is discrete 
and, since $\bar\Ga^o$ is not nilpotent, 
we have $\dim Z_k<\dim\bar\Ga$. 
Proposition \ref{growth}
now follows by applying the next lemma 
with $H=\bar\Ga$ and $H_1=Z_k$. 
$\Box$

\begin{lem}
Let $H$ be a Lie group, let $H_1\lhd H$ be a closed normal
subgroup so that $\bar H:=H/H_1$ is a positive dimensional
Lie group with discrete center, and suppose $\Ga\subset H$
is a dense, finitely generated subgroup.  If $U$ is 
any neighborhood of $e$ in $H$, then the function
$f(k):=\#\{ g\in \Ga : \mbox{$|g|_\Ga\leq k$, $g\in U$}\}$
grows superpolynomially.
\end{lem}

\proof
The idea of the proof is to use the contracting property of 
commutators to produce a sequence $\{\al_k\}$
in $H\cap \Ga$ which converges exponentially to the identity.
The word norm $|\al_k|_\Ga$ grows exponentially with $k$, 
but the number of elements of $\langle\al_1,\ldots,\al_k\rangle$ in 
$U$ also grows exponentially with $k$; by comparing growth
exponents we find that $f$ grows superpolynomially.

Fix $M\in\N$, a positive real number $\eps<1/3$  
and some left-invariant Riemannian metric on $H$. 
Since the differential of the commutator map 
$(h,h')\mapsto[h,h']$ vanishes at $(e,e)$ 
we can find a neighborhood $V$ of $e$ in $H$
such that:
\begin{equation}
\label{eny}
h,h'\in V \quad\Lra\quad [h,h']\in V \quad\hbox{and}\quad 
d([h,h'],e)<\frac{1}{2M} d(h,e) 
\end{equation}
Since the differential of the $k$-th power $h\mapsto h^k$ at $e$ 
is $k\cdot id_{T_eH}$ for all $k\in\Z$, 
we can furthermore acheive that, 
whenever $1\leq k,k'\leq M$ and $h,h^k,h^{k'}\in V$,  
then 
\begin{equation}
\label{oqa}
d(h^k,h^{k'}) \geq (|k-k'|-\eps)\cdot d(h,e)
\end{equation}
By our assumption, there exist finitely many elements 
$\ga_1,\dots,\ga_m\in\Ga\cap V$ such that the centralizers 
$Z_{\bar H}(\bar\ga_j)$ of their images in $\bar H$ have discrete intersection. 
We construct an infinite sequence of elements 
$\al_i\in(\Ga\cap V)\setminus H_1$ by picking $\al_0\in V$ arbitrarily 
and setting 
$\al_{i+1}=[\al_i,h_{j(i)}]\not\in H_1$ for suitably chosen $1\leq j(i)\leq m$. 
Then
\begin{equation}
\label{pwv}
0<d(\al_{i+1},e)<\frac{1}{2M}d(\al_i,e)
\end{equation}
by (\ref{eny}).
\begin{sublem}
Pick $n_0\in\N$. 
The $M^n$ elements 
\begin{equation}
\label{theelements}
\ga_{\eps_1\dots\eps_n}=\al_{n_0+1}^{\eps_1}\cdots\al_{n_0+n}^{\eps_n}
\qquad \eps_i\in\{0,\dots,M-1\}
\end{equation}
are distinct.
\end{sublem}
\proof
Assume that $\ga_{\eps_1\dots\eps_n}=\ga_{\eps'_1\dots\eps'_n}$,
$\eps_l\neq\eps'_l$ and $\eps_i=\eps'_i$
for all $i<l$. 
Then
\[
\al_{n_0+l}^{\eps_l-\eps'_l}=
\al_{n_0+l+1}^{\eps'_{l+1}-\eps_{l+1}}\cdots\al_{n_0+n}^{\eps'_n-\eps_n}.
\]
On the other hand (\ref{oqa},\ref{pwv}) 
and the triangle inequality imply 
\[
d(\al_{n_0+l+1}^{\eps'_{l+1}-\eps_{l+1}}\cdots\al_{n_0+n}^{\eps'_n-\eps_n},e)
 < M\cdot \sum_{j=1}^{\infty}\frac{1}{(2M)^j} \cdot d(\al_{n_0+l},e) 
 < \frac{1}{2} d(\al_{n_0+l},e) 
 < d(\al_{n_0+l}^{\eps_l-\eps'_l},e) ,
\]
a contradiction. 
\qed 

To complete the proof of the lemma, 
we observe that 
the elements (\ref{theelements}) have word norm 
$|\ga_{\eps_1\dots\eps_n}|_{\Ga}\leq const(n_0)\cdot 2^n$ 
and are contained in $U$ if $n_0$ is sufficiently large. 
This shows that $f(k)$ grows polynomially of order at least
$\frac{log(M)}{log(2)}$ for all $M$,
hence the claim.
\qed

\section{Proof of the main theorem}

\no
{\em Proof of theorem \ref{mainthm}:}
Let $\rho_0:\Ga\times\Ga\ra\Ga$ be the isometric action
of $\Ga$ on itself by left translation, and 
let $\phi:\Ga\ra Nil\times X$ be a quasi-isometry.
Then there is a quasi-action $\rho$ of $\Ga$ on $Nil\times X$
such that $\phi$ quasi-isometrically conjugates $\rho_0$ into $\rho$. 
According to section \ref{projectingqactions},
$\rho$ projects to a cobounded quasi-action $\bar\rho$ of $\Ga$ on $X$. 
After quasi-isometrically conjugating $\rho$ we may assume that $\bar\rho$ 
is an isometric action, cf.\ \ref{conjqactions}. 
For any point $y\in Nil\times X$ the map 
$\Ga\ra Nil\times X$ given by $\ga\mapsto\rho(\ga, y)$ is a quasi-isometry. 
Thus for all $x\in X$ and $R>0$
$$\#\{ \ga\in\Ga\mid 
\mbox{$|\ga|_\Ga < k$, $\bar\rho(\ga)\cdot x\in B_R(x)$}\}$$
grows at most as fast as the volume of balls in $Nil$,
i.e. $<k^d$ for some $d$. 
Proposition \ref{growth} implies that 
$L:=\bar\rho(\Ga)$ is a discrete subgroup in $Isom(X)$ 
and hence a uniform lattice. 
The kernel $H$ of the action $\bar\rho$ is then quasi-isometric 
to the fiber $Nil$. 
$\Box$

\noindent

\medskip\no 
Bruce Kleiner\\
Department of Mathematics\\
University of Pennsylvania\\ 
Philadelphia, PA 19104-6395\\ 
bkleiner$@$math.upenn.edu

\medskip\no
Bernhard Leeb\\
Mathematisches Institut\\
Universit\"at Bonn, Beringstr.\ 1\\
53115 Bonn, Germany\\
leeb$@$rhein.iam.uni-bonn.de

\end{document}